\documentstyle[12pt,leqno,amssymb,graphics]{amsart}

\newtheorem{thm}{Theorem}[section]
\newtheorem{lemma}[thm]{Lemma}
\newtheorem{prop}[thm]{Proposition}
\newtheorem{cor}[thm]{Corollary}

\theoremstyle{definition}
\newtheorem{dfn}[thm]{Definition}
\theoremstyle{remark}

\begin{document}

\newcommand{\ct}{\cite}
\newcommand{\pr}{\protect\ref}
\newcommand{\su}{\subseteq}
\newcommand{\pa}{{\partial}}
\newcommand{\im}{{Imm(F,\E)}}
\newcommand{\lm}{{\lambda}}
\newcommand{\tc}{{\mathfrak{T}}}
\newcommand{\Q}{{\Bbb Q}}
\newcommand{\hf}{{1 \over 2}}
\newcommand{\CC}{{\mathcal C}}

\newcommand{\R}{{\Bbb R}}
\newcommand{\Z}{{\Bbb Z}}
\newcommand{\E}{{{\Bbb R}^3}}
\newcommand{\C}{{{\Bbb Z}/2}}
\newcommand{\B}{{\Bbb B}}

\newcommand{\FF}{{\mathcal F}}

\newcommand{\A}{{\mathcal{A}}}
\newcommand{\Ai}{{{\mathrm{Arf}}(g^i)}}
\newcommand{\Ar}{{\mathrm{Arf}}}
\newcommand{\r}{{\mathrm{rank}}}
\newcommand{\I}{{\mathrm{Id}}}
\newcommand{\spp}{{\mathrm{supp}}}

\newcommand{\F}{{\mathbf{F}}}
\newcommand{\G}{{\Bbb G}}
\newcommand{\cd}{{\mathrm{codim}}}
\newcommand{\p}{{\psi}}
\newcommand{\hp}{{\Psi}}
\newcommand{\hhp}{\widehat{\Psi}}
\newcommand{\pp}{{\widehat{\psi}}}

\newcommand{\4}{{\mathcal{H}}}
\newcommand{\N}{{\mathcal{N}}}
\newcommand{\ng}{{{\mathcal{N}}_g}}

\newcommand{\hM}{{\widehat{\mathcal{M}}}}
\newcommand{\m}{{\mathcal{M}}_g}
\newcommand{\hm}{{\hM_g}}
\newcommand{\hmi}{{\hM_{g^i}}}
\newcommand{\mi}{{\M_{g^i}}}

\newcommand{\f}{{\overline{f}}}

\newcommand{\hc}{{H_1(F;\C)}}
\newcommand{\ohg}{{O(\hc,g)}}

\newcommand{\ohi}{{O(\hc,g^i)}}
\newcommand{\Tc}{{{\mathcal{T}}_c}}
\newcommand{\SP}{{{\mathcal{S}}_P}}

\newcommand{\Y}{{\mathcal{Y}}}
\newcommand{\T}{{\mathcal{T}}}
\newcommand{\U}{{\mathcal{U}}}
\newcommand{\tF}{{\widetilde{F}}}

\newcommand{\one}{{\widetilde{F_1}}}
\newcommand{\two}{{\widetilde{F_2}}}
\newcommand{\tfk}{{\widetilde{F_k}}}
\newcommand{\th}{{\widetilde{h}}}
\newcommand{\ep}{{\epsilon}}
\newcommand{\wb}{{W^\bot}}
\newcommand{\w}{{\bigcup_i w_i}}
\newcommand{\eep}{{\widehat{\epsilon}}}

\newcounter{numb}

\title{Immersions of Non-orientable Surfaces}
\author{Tahl Nowik}
\address{Department of Mathematics, Bar-Ilan University, 
Ramat-Gan 52900, Israel.}
\email{tahl@@math.biu.ac.il}
\date{October 1, 2003}
\thanks{Partially supported by the Minerva Foundation}

\begin{abstract}
Let $F$ be a closed non-orientable surface.
We classify all finite order invariants of immersions of $F$
into $\E$, with values in any Abelian group. We show they are all
functions of the universal order 1 invariant that we construct 
as $T \oplus P \oplus Q$ where $T$ is a $\Z$ valued invariant reflecting the number of triple
points of the immersion, and $P,Q$ are $\C$ valued invariants characterized by the property
that for any regularly homotopic 
immersions $i,j:F\to \E$, $P(i)-P(j) \in \C$ 
(respectively $Q(i)-Q(j) \in \C$)
is the number mod 2 of tangency points (respectively quadruple points)
occurring in any generic regular homotopy between $i$ and $j$. 

For immersion $i:F\to\E$ and diffeomorphism $h:F\to F$ such that $i$ and $i \circ h$
are regularly homotopic we show:
$$P(i\circ h)-P(i)=Q(i\circ h)-Q(i) =\bigg(\r(h_*-\I)+\ep(\det h_{**})\bigg)\bmod{2}$$
where $h_*$ is the map induced by $h$ on $\hc$, $h_{**}$ is
the map induced by $h$ on $H_1(F;\Q)$, and for $0 \neq q \in \Q$, 
$\ep(q) \in\C$ is 0 or 1
according to whether $q$ is positive or negative, respectively.
\end{abstract}

\maketitle

\section{Introduction}

Finite order invariants of immersions of a closed surface into $\E$ have been defined in \ct{o},
and the case of orientable surfaces has been studied in \ct{a},\ct{o},\ct{h},\ct{f}.
In the present work we establish all analogous results for the non-orientable case.
We classify all finite order invariants. For each $n$ we construct a
universal order $n$ invariant, and for any $n>1$ it
is constructed as an explicit function of the universal order 1 invariant.
The universal order 1 invariant is given by 
$T \oplus P \oplus Q$ where $T$ is a $\Z$ 
valued invariant reflecting the number of triple
points of the immersion, and $P,Q$ are $\C$ valued invariants characterized by the property
that for any regularly homotopic 
immersions $i,j:F\to \E$, $P(i)-P(j) \in \C$ 
(respectively $Q(i)-Q(j) \in \C$)
is the number mod 2 of tangency points (respectively quadruple points)
occurring in any generic regular homotopy between $i$ and $j$. 
We give an explicit formula for $P(i \circ h)-P(i)$ and $Q(i \circ h)-Q(i)$ 
for any diffeomorphism $h:F \to F$ such that $i$ and $i \circ h$ are regularly homotopic.
This requires the study of the quadratic form induced on $\hc$ by an immersion $i:F\to\E$,
where in the case of non-orientable surfaces we name it an $\4$-form, rather than quadratic form,
to emphasize its distinct character.  In order to obtain the stated formulae for $P,Q$, we study the 
group of diffeomorphisms of $F$ which preserves the induced $\4$-form.

The extension of the classification of finite order invariants from the orientable to the 
non-orientable case (Section \pr{fin}) is rather straight forward, whereas
the study of the group of diffeomorphisms preserving an $\4$-form (Section \pr{og}), and the formula 
for the invariants $P,Q$ (Section \pr{fmq}), will require substantially new analysis.

\section{Finite order invariants}\label{fin}

Let $F$ be a closed non-orientable surface. $\im$ denotes the 
space of all immersions of $F$ into $\E$, with the $C^1$ topology.
A CE point of an immersion $i:F\to \E$ is a point of self intersection
of $i$ for which the local stratum in $\im$ corresponding to the 
self intersection, has codimension one. 
For $F$ non-orientable we distinguish four types of CEs which we name
$E, H, T, Q$. 
In the notation of \ct{hk} they are respectively
$A_0^2|A_1^+$, $A_0^2|A_1^-$, $A_0^3|A_1$, $A_0^4$.
The four types may be demonstrated by the following local models, where letting $\lm$ vary, 
we obtain a 1-parameter family of immersions which is transverse to the 
given codim 1 stratum, intersecting it at $\lm=0$.

$E$: \ \ $z=0$, \ \ $z=x^2+y^2+\lm$. See Figure \pr{fet}, ignoring the vertical plane.

$H$: \ \ $z=0$, \ \ $z=x^2-y^2+\lm$. See Figure \pr{fht}, ignoring the vertical plane.

$T$: \ \ $z=0$, \ \ $y=0$, \ \ $z=y+x^2+\lm$. See Figure \pr{ftq}, ignoring the vertical 
plane $x=0$.

$Q$: \ \ $z=0$, \ \ $y=0$, \ \ $x=0$, \ \ $z=x+y+\lm$. This is simply four planes passing 
through one point, any three of which are in general position.

A choice of one of the two sides of the local codim 1 stratum at a given point of the 
stratum, is represented by the choice of $\lm<0$ or $\lm>0$ in the formulae above. We will 
refer to such a choice as a 
\emph{co-orientation} for the configuration of the self intersection. 
For types $E$ and $T$, the configuration of the self intersection at the two sides 
of the stratum is distinct, namely, for $\lm<0$ there is an additional 
2-sphere in the image of the immersion, and we permanently
choose this side ($\lm<0$) as our positive side for the co-orientation. 
For types $H$ and $Q$, the configuration of the self intersection on the two 
sides of the strata are 
indistinguishable, and in fact we will see that the strata in this case are globally one-sided,
and so no coherent choice of co-orientation is possible.

We fix a closed non-orientable surface $F$ and a regular homotopy class $\A$ of 
immersions of $F$ into $\E$.
We denote by $I_n\su \A$ ($n\geq 0$) the space of all immersions in $\A$ which have precisely 
$n$ CE points (the self intersection being elsewhere stable).
In particular, $I_0$ is the space of all stable immersions in $\A$. 

For an immersion $i:F\to\E$ having a CE located at $p\in\E$, we will now define
the degree $d_p(i)$. This will differ from the definition
for orientable surfaces, given in \ct{o}, in two ways. 
First, for nonorientable surface, 
the degree of a map $F \to \E - \{ p\}$ is defined only mod 2. Second,
we do not have an orientation which determines into what side of $p$ we must 
push each sheet participating in the CE,
in order for the degree to be computed. 
So for non-orientable surface we define $d_p(i) \in \{ + , - \}$, as follows: 
For CE of type $E,T$ we move the immersion to the positive side determined by its
permanent co-orientation, obtaining an immersion $i'$ where a new 
2-sphere appears in the image. We define $d_p(i)\in \{ + , - \}$ to be the degree of the 
map $i':F \to \E- \{ p' \}$ where $p'$ is any point in the open 3-cell bounded by the new 2-sphere,
and where $+ =$ even and $- =$ odd. 
For CEs of type $H,Q$,
let $V$ denote the region between the sheets of the surface which appears once we move 
away from the stratum (i.e. once $\lm\neq 0$). 
For type $Q$ this is a 3-simplex defined by our four sheets. For type $H$ this region is not 
bounded by the local configuration,
but may still be defined, e.g. for $\lm>0$ in the formula for $H$, 
$V$ will be a region consisting of points close to the origin and 
satisfying $0\leq z\leq x^2-y^2+\lm$. We define 
$d_p(i)\in \{ + , - \}$ to be the degree of the 
map $i':F \to \E- \{p' \}$ where $p'$ is any point in $V$, and notice that since there is an even
number of sheets involved in a CE of type $H$ or $Q$, the definition is independent of the side of the 
stratum we choose to move into, (i.e. whether we choose $\lm>0$ or $\lm<0$).

We define
$C_p(i)$ to be the expression $R_e$ where $R$
is the symbol describing the configuration of the CE of $i$ at $p$ ($E$, $H$, $T$, or $Q$) and $e=d_p(i)$. 
We define $\CC_n$ to be the set of all
\emph{un-ordered} $n$-tuples of expressions $R_e$ with $R$ one of the four symbols
and $e\in\{+,-\}$.
So $\CC_n$ is the set of un-ordered $n$-tuples of elements of 
$\CC_1 = \{E_+ , E_- , H_+ , H_- , T_+ , T_- , Q_+ , Q_- \}$.
We define $C:I_n \to \CC_n$ by
$C(i)=[C_{p_1}(i),\dots, C_{p_n}(i)] \in \CC_n$ where $p_1,\dots,p_n$ are the $n$
CE points of $i$. The map 
$C:I_n \to \CC_n$ is easily seen to be surjective.

A regular homotopy between two immersions in $I_n$ is called an AB equivalence if 
it is alternatingly of type A and B, where $J_t:F\to\E$ ($0\leq t \leq 1$)
is of type A if it is of the form
$J_t = U_t \circ i \circ V_t$ where 
$i:F\to \E$ is an immersion and
$U_t:\E\to\E$, 
$V_t:F\to F$ are isotopies, and $J_t:F\to\E$ ($0\leq t \leq 1$)
is of type B if $J_0 \in I_n$ and there are  
little balls $B_1,\dots,B_n\su\E$ centered at the $n$ CE points of $J_0$
such that $J_t$ fixes $U=(J_0)^{-1}(\bigcup_k B_k)$ 
and moves $F-U$ within $\E - \bigcup_k B_k$.

\begin{prop}\label{AB}
Let  $i,j\in I_n$, then $i$ and $j$ are AB equivalent iff $C(i)=C(j)$.
\end{prop}

\begin{pf}
The proof proceeds as in \ct{o} Proposition 3.4 except for one step that must be added 
in the present case, of non-orientable surfaces. 
In the final stage of the proof of \ct{o} Proposition 3.4, 
we have that $i'''|_D$ and $j|_D$ 
are homotopic in $\E-\bigcup_k B_k$ relative $\pa D$, (since 
$\widehat{d}_{p_k}(i''')=\widehat{d}_{p_k}(j)$ for all $k=1,\dots,n$,
where $\widehat{d}_p$ denotes the $\Z$ valued degrees defined in \ct{o} for the orientable case).
In the present case this will not necessarily be true. If $f:S^2 \to \E-\bigcup_k B_k$ 
is the map determined by the pair of maps $i'''|_D,j|_D$, 
then instead of $f$ being null-homotopic 
as in the orientable case, it may only be 
of even degree with respect to each $p_1,\dots,p_n$.
This can be remedied as follows: Let $h$ be one of the 1-handles having both
ends attached to $D_1$, so
$D_1 \cup h$ is a Mobius band and the 2-handle $D$ is
glued to it twice in the same direction.
A homotopy of $h$ which traces a sphere in $\E$ enclosing just one $p_i$ and not the others, will
change the degree of $f$ with respect to $p_i$ by 2, and leave the degree with respect
to all other $p_k$ unchanged. We realize such homotopies by regular homotopies
until the new map $f$ has degree 0 with respect to each $p_1,\dots,p_n$.

\end{pf}

Given an immersion $i\in I_n$, a \emph{temporary co-orientation} for $i$ is a choice 
of co-orientation at each of the $n$ CE points $p_1, \dots , p_n$ of $i$.
Given a temporary co-orientation $\tc$ for $i$ and a subset $A\su \{p_1,\dots,p_n\}$,
we define $i_{\tc,A} \in I_0$ to be the immersion obtained from $i$ by resolving all CEs
of $i$ at points of $A$ into the 
positive side with respect to $\tc$,
and all CEs not in $A$ into the negative side.
Now let $\G$ be any Abelian group and let $f:I_0\to\G$ be an invariant, i.e. a function which is 
constant on each connected component of $I_0$.
Given $i\in I_n$ and a temporary co-orientation $\tc$ for $i$,
$f^\tc(i)$ is defined as follows:
$$f^\tc(i)=\sum_{ A \su \{p_1,\dots,p_n\} } (-1)^{n-|A|} f(i_{\tc,A})$$
where $|A|$ is the number of elements in $A$.
The statement $f^\tc(i)=0$ is independent of the temporary co-orientation $\tc$
so we
simply write $f(i)=0$.
An invariant $f:I_0\to\G$ is called \emph{of finite order} if 
there is an $n$ such that $f(i)=0$ for all $i\in I_{n+1}$.
The minimal such $n$ is called the \emph{order} of $f$.
The group of all invariants on $I_0$ of order at most $n$ is denoted 
$V_n = V_n(\G)$.

Let $f \in V_n$. If $i\in I_n$ has at least one CE
of type $H$ or $Q$ and $\tc$ is a temporary co-orientation for $i$,
then $2f^\tc(i)=0$, the proof being the same as in \ct{o} Proposition 3.5.
So in this case $f^\tc(i)$ is independent of $\tc$. 
We use this fact to extend any $f\in V_n$ to $I_n$
by setting for any $i\in I_n$, $f(i) = f^\tc(i)$,
where if $i$ includes at least one CE of type $H$ or $Q$ then
$\tc$ is arbitrary, and if all CEs of $i$ are of type $E$ and $T$
then the permanent co-orientation is used for all CEs of $i$.
We will always assume without mention that any $f \in V_n$ is extended to $I_n$ 
in this way.
(If $f\in V_n$ then we are not extending $f$ to $I_k$ for
$0<k<n$). 

We remark at this point that the same argument as explained in \ct{o} Remark 3.7, showing that
for orientable surfaces, the strata corresponding to configurations $H^1$ and $Q^2$
may not be globally co-oriented, will show that
the same is true for non-orientable surfaces for configurations $H$ and $Q$.

For $f\in V_n$ and $i,j\in I_n$, 
if $C(i)=C(j)$ then $f(i)=f(j)$, the proof being the same as in \ct{o} Proposition 3.8,
so any $f\in V_n$ induces a well defined
function $u(f):\CC_n\to\G$. 
The map $f\mapsto u(f)$ induces an injection $u:V_n / V_{n-1} \to \CC_n^*$ 
where $\CC_n^*$ is the group of all functions from $\CC_n$ to $\G$.
We will find the image of $u$ for all $n$, by this we classify all finite
order invariants.

Let $i\in \A$ be an immersion with a self intersection of local 
codim 2 at $p$ and 
$n-1$ additional self-intersections of local codim 1 (i.e. CEs) at $p_1,\dots,p_{n-1}$.
We look at a 2-parameter family of immersions which moves $F$ only in a neighborhood 
of $p$, such that the immersion $i$ corresponds to 
parameters $(0,0)$ and such that this 2-parameter family is transverse to the local 
codim 2 stratum at $i$.
In this 2-parameter family of immersions we look at a loop which 
encircles the point of 
intersection with the codim 2 strata,
i.e. a circle around the origin in the parameter plane.
This circle crosses the local codim 1 strata some $r$ times.
Between each two intersections we have an immersion in $I_{n-1}$ 
with the same $n-1$ CEs, at $p_1,\dots,p_{n-1}$.
At each intersection with the local codim 1 strata, an $n$th CE is added, 
obtaining an immersion in $I_n$. Let $i_1,\dots,i_r$ be the $r$ immersions in $I_n$ so 
obtained and let $\ep_k$, $k=1,\dots,r$ be $1$ or $-1$ according to whether we are passing 
the $n$th CE of $i_k$ in the direction of its 
permanent co-orientation, if it has one, and if the CE is of type $H$ or $Q$ then $\ep_k$ is 
arbitrarily chosen. 
For $f \in V_n$, $f(i_k)$ is defined and it is easy to see that 
$\sum_{k=1}^r \ep_k f(i_k) = 0$. 
Looking at $u:V_n / V_{n-1} \to \CC_n^*$  we thus obtain relations that must be satisfied by 
a function in $\CC_n^*$ in order for it to lie in the image of $u$.
We will now find all relations on $\CC_n^*$ obtained in this way.
The relations on a $g\in\CC_n^*$ will be written as relations on the symbols $R_e$, e.g.
$0 = T_e - T_{-e}$ will stand for the set of all relations of the form
$0 = g([T_e , {R_2}_{e_2},\dots,{R_n}_{e_n}]) - 
g([T_{-e}, {R_2}_{e_2},\dots,{R_n}_{e_n}])$ with 
arbitrary ${R_2}_{e_2},\dots,{R_n}_{e_n}$. 
We already know that the following two relations hold:
$0=2H_e$, $0=2Q_e$ (for both $e$).

We look at local 2-parameter 
families of immersions which are transverse to 
the various local codim 2 strata.
These may be divided into six types which we name after the types of CEs 
appearing in the 2-parameter family:
$EH$, $TT$, $ET$, $HT$, $TQ$, $QQ$. In the notation of \ct{hk} they are 
respectively: $A_0^2|A_2$, $A_0^3|A_2$, $(A_0^2|A_1^+)(A_0)$, 
$(A_0^2|A_1^-)(A_0)$, $(A_0^3|A_1)(A_0)$, $A_0^5$.
Formula and sketch for local model for such strata,
the bifurcation diagrams, and the relations obtained, are as follows.
A sign $\pm$ appears wherever the element is known to be of order 2 (which is
wherever there is no co-orientation).

\begin{figure}[h]
\scalebox{0.6}{\includegraphics{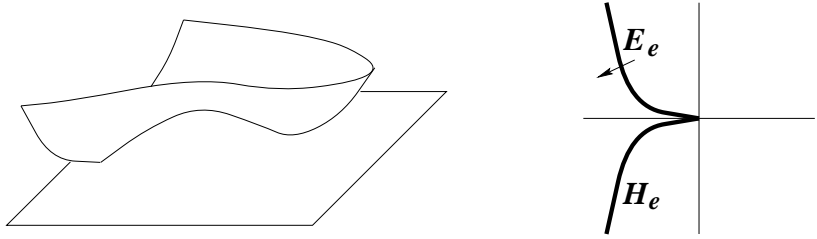}}
\caption{$EH$ configuration}\label{feh}
\end{figure}

$EH$: \ \ $z=0$, \ \ $z=y^2 + x^3+\lm_1 x + \lm_2$.
\begin{equation}\label{eeh}
0 = E_e \pm H_e
\end{equation}

\begin{figure}[h]
\scalebox{0.6}{\includegraphics{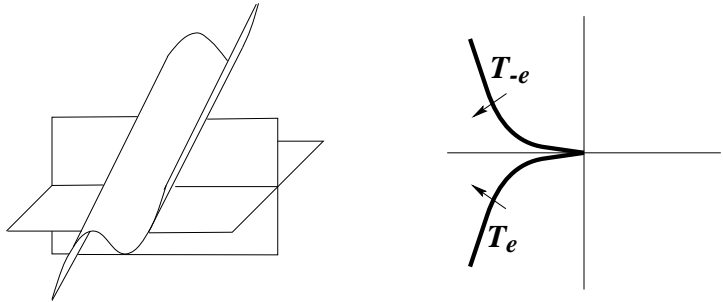}}
\caption{$TT$ configuration}\label{ftt}
\end{figure}

$TT$: \ \ $z=0$, \ \ $y=0$, \ \ $z=y+x^3+\lm_1 x  + \lm_2$.
\begin{equation}\label{ett}
0 = T_e - T_{-e}
\end{equation}

\begin{figure}[h]
\scalebox{0.6}{\includegraphics{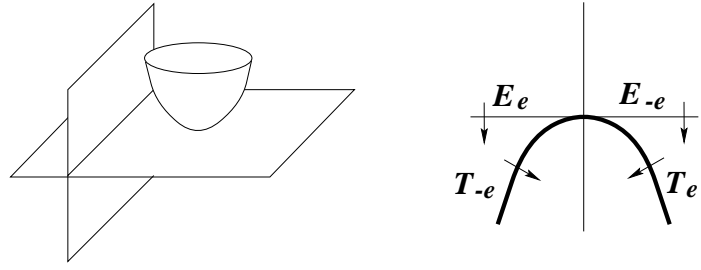}}
\caption{$ET$ configuration}\label{fet}
\end{figure}

$ET$: \ \ $z=0$, \ \ $x=0$, \ \ $z=(x-\lm_1)^2 + y^2 + \lm_2$.
\begin{equation}\label{eet}
0 = T_{-e} - T_e - E_{-e} + E_e
\end{equation}

\begin{figure}[h]
\scalebox{0.6}{\includegraphics{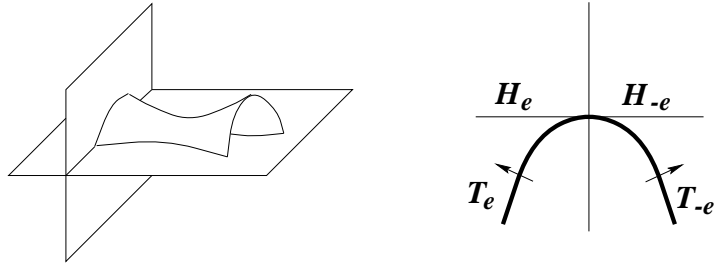}}
\caption{$HT$ configuration}\label{fht}
\end{figure}

$HT$: \ \ $z=0$, \ \ $x=0$, \ \ $z=(x-\lm_1)^2 - y^2 + \lm_2$.
\begin{equation}\label{eht}
0 = -T_e + T_{-e} \pm H_{-e} \pm H_e
\end{equation}

\begin{figure}[h]
\scalebox{0.6}{\includegraphics{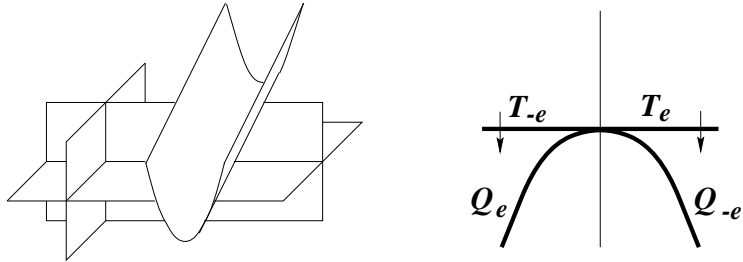}}
\caption{$TQ$ configuration}\label{ftq}
\end{figure}

$TQ$: \ \ $z=0$, \ \ $y=0$, \ \ $x=0$, \ \ $z=y+(x-\lm_1)^2 + \lm_2$.
\begin{equation}\label{etq}
0 =  \pm Q_e \pm Q_{-e} - T_e + T_{-e}
\end{equation}

\begin{figure}[h]
\scalebox{0.6}{\includegraphics{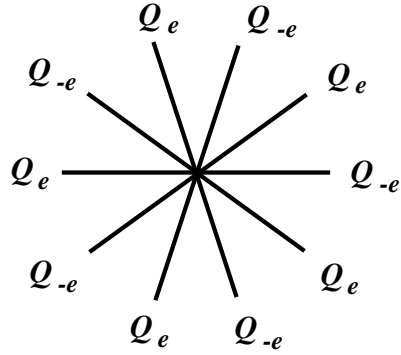}}
\caption{$QQ$ configuration}\label{fqq}
\end{figure}

$QQ$: \ \ Five planes meeting at a point.
\begin{equation}\label{eqq}
0 = 5(\pm Q_e) + 5(\pm Q_{-e})
\end{equation}

For the first five types,
the bifurcation diagram and degrees are obtained from the sketch and formula in a straight forward 
manner. 
The diagram for $QQ$ is obtained as explained in \ct{o}.
The degrees appearing in the present $QQ$ diagram require explanation, but we may avoid it since 
from the first five relations it already follows that $Q_+=Q_-$ and so since $2Q_e=0$,
the last relation does not add to the first five relations, whatever the degrees may be. 
Letting $\B = \{ x \in \G : 2x=0 \}$ the above relations may be summed up as follows:

\begin{itemize}
\item $T_+ = T_-$
\item $E_+ = E_- = H_+ = H_- \in \B$
\item $Q_+ = Q_- \in \B$.
\end{itemize}
We denote by $\Delta_n = \Delta_n(\G)$ the subgroup of $\CC_n^*$ of all functions 
satisfying these relations, so the image of $u:V_n \to \CC_n^*$ is contained in 
$\Delta_n$.

We define the universal Abelian group $\G_U$ by the Abelian group presentation
$\G_U = \left< t , p , q \ | \ 2p = 2q = 0 \right>$. 
We define the universal element $g_1^U\in\Delta_1(\G_U)$ by $g_1^U(T_e) = t$,
$g_1^U(E_e)=g_1^U(H_e)=p$, $g_1^U(Q_e)=q$, 
so indeed $g_1^U\in\Delta_1(\G_U)$. 
Then for arbitrary Abelian group $\G$ we have 
$\Delta_1(\G) \cong Hom(\G_U , \G)$ where the isomorphism 
maps a homomorphism $\phi\in Hom(\G_U , \G)$ 
to the function $\phi \circ g_1^U \in \Delta_1(\G)$.
We will show that there is an order 1 
invariant $f_1^U:I_0\to\G_U$ with $u(f_1^U)=g_1^U$. It will follow that for any 
group $\G$, $u:V_1/V_0\to \Delta_1(\G)$ is surjective, since if $g\in\Delta_1(\G)$ 
and $g=\phi \circ g_1^U$ where $\phi\in Hom(\G_U , \G)$ 
then $u(\phi\circ f_1^U) = g$.

We define $f_1^U:I_0\to\G_U$ as follows:
Choose a base immersion $i_0 \in I_0$ once and for all.
Then for each $i \in I_0$ take a generic regular homotopy
$H_t:F \to \E$ ($0 \leq t \leq 1$)
from $i_0$ to $i$ and define $f_1^U(i)=n_1 t + n_2 p + n_3 q \in \G_U$ where
$n_1 \in \Z$ is the number of CEs of type $T$ occurring along $H_t$, each counted as $\pm 1$ according
to the permanent coorientation of $T$, and $n_2 \in \C$ (respectively $n_3 \in \C$) is the number mod 2 of
tangencies (respectively quadruple points), occurring in $H_t$.
If $f_1^U$ is well defined, i.e. independent of the choice of $H_t$, then clearly $u(f_1^U)=g_1^U$.
So it remains to show that $f_1^U$ is independent of $H_t$, which is equivalent to showing that
the value is 0 for any \emph{closed} regular homotopy, i.e. one that begins and ends with the 
same immersion.

Now $\pi_1(\A) \cong \C \oplus \C$. 
Though this is different from the fundamental group in the orientable case ($\C \oplus \Z$),
the two generators of $\pi_1(\A)$ are obtained in the same way, namely, the first generator is given by
one full rigid rotation of the surface in $\E$, and the second generator is obtained by
only moving a disc in $F$ using the infinite cyclic generator in $\pi_1(Imm(S^2,\E))$
as shown in \ct{q}. So the proof that $f_1^U$ is well defined may proceed exactly as in the orientable case,
as appears in \ct{o}, once we know that the same invariant is well defined for $S^2$.
But for $S^2$, \ct{o} gives a full classification, and if hats will denote the corresponding constructions
for orientable surfaces, then this particular invariant for $S^2$ is obtained from the following
$g \in \widehat{\Delta}_1(\G_U)$ (note we are taking orientable $\Delta_1$ of non-orientable $\G_U$):
$g(T^a_m)=t$, $g(E^a_m)=g(H^a_m)=p$, $g(Q^a_m)=q$, for all $a,m$. (See \ct{o} for definition of the symbols
$R^a_m$).
This completes the proof that $f_1^U$ is well defined.
The invariant $f_1^U$ is clearly a \emph{universal} order 1 invariant, 
where universal order $n$ invariant is defined as follows:

\begin{dfn}\label{uni}
A pair $(\G,f)$ where $\G$ is an Abelian group and 
$f:I_0 \to \G$ is an order $n$ invariant,
will be called a \emph{universal order $n$ invariant}
if for any Abelian group $\G'$ and any order $n$ invariant 
$f':I_0 \to \G'$ there exists a unique homomorphism $\varphi:\G \to \G'$ 
such that $f' - \varphi \circ f$ is an invariant of order 
at most $n-1$. 
\end{dfn}

We have $\G_U  = \Z t \oplus (\C)p \oplus (\C)q$ and let $T,P,Q$
be the projection of $f_1^U$ into these three factors of $G_U$. 
We may write an explicit formula for $T$, namely, $T(i) = {N-c \over 2}$ 
where $N$ is the number of triple points in $i$
and $c=0$ if $\chi(F)$ is even, $c=1$ if $\chi(F)$ is odd.
Indeed, passing a CE of type $T$ in the positive direction
with respect to the permanent coorientation increases the number of triple point by 2, and the other CE
types leave the number of triple points unchanged.
(This formula may however differ by a constant from the $T$ defined using the given $i_0$ as base immersion).
Note that indeed $N-c \over 2$ is an integer, as is shown in 
\ct{b}, and may also be deduced using our present considerations.
For $P,Q$ we do not have an explicit formula for all immersions, however in Section \pr{fmq}
we will give an explicit formula for a certain situation.

We now classify all higher order invariants. We will define $E_n \su \Delta_n$ by two additional
restrictions on the functions $g\in \Delta_n$.
Let $Y = \{ T_+,H_+,Q_+ \} \su \CC_1$, then any $g \in \Delta_n$ is determined by its values on
un-ordered $n$-tuples of elements of $Y$ and so we may state these relations in terms of
such $n$-tuples. 
Given an un-ordered $n$-tuple $z$ of elements of $Y$, 
we define $m_{H_+}(z)$ and $m_{Q_+}(z)$ 
as the number of times that
$H_+$ and $Q_+$ appear in $z$ respectively.
We define $r(z)$, (the \emph{repetition} of $H_+$ and $Q_+$ in $z$), as 
$$r(z)= \max(0, m_{H_+}(z)-1) + \max(0, m_{Q_+}(z)-1).$$

\begin{dfn}\label{En}
Given an Abelian group $\G$,
$E_n = E_n(\G) \su \Delta_n(\G)$ is the subgroup consisting
of all $g\in \Delta_n(\G)$
satisfying the following two additional restrictions:

\begin{enumerate}

\item When $n \geq 3$, $g$ must satisfy the relation $H_+ H_+ Q_+ = H_+ Q_+ Q_+$.
By this we mean that 
$g([H_+, H_+, Q_+, {R_4}_{e_4},\dots,{R_n}_{e_n}]) =
g([H_+, Q_+, Q_+, {R_4}_{e_4},\dots,{R_n}_{e_n}])$ for
arbitrary ${R_4}_{e_4},\dots,{R_n}_{e_n} \in Y$.
\label{e1}

\item For any un-ordered $n$-tuple $z$ of elements of $Y$,
$g(z)\in 2^{r(z)}\G$, i.e. there exists an element $a\in \G$ such
that $g(z)=2^{r(z)} a$.
\label{e2}

\end{enumerate}
\end{dfn}

We define algebraic structures $K \su L \su M$, where $L$ is a commutative 
ring, $K$ is a subring of $L$, and $M$ is a module over $K$.
$L$ is defined as the ring of formal power series with integer coefficients
and variables $t,p,q$ and with relations 
\begin{itemize}
\item $2p=2q=0$.
\item $p^2 q = p q^2$. 
\end{itemize}

Given a monomial $f$, 
we define $m_p(f)$ and $m_q(f)$ as the multiplicity of $p$ and $q$ in $f$ respectively. We 
define $r(f)$, (the repetition of $p$ and $q$ in $f$), as 
$$r(f)= \max(0, m_p(f)-1) + \max(0, m_q(f)-1).$$ 
$r(f)$ is preserved under the 
relations in $L$ and so is well defined on 
equivalence classes of monomials.
The equivalence class of a monomial $f$ 
will be denoted $\f$. 
Now $K \su L$ is defined to be the subring of power series
including only the variable $t$.
On the other hand 
we extend $L$ to a larger structure $M$ which will 
be a module over the subring $K$, 
as follows:
For each $\f \in L$ for which $f$ is a monomial with coefficient 1,
we adjoin a new element $\zeta_\f$ satisfying the
relation $2^{r(\f)} \zeta_\f = \f$.  
Now we let $K$ act on $M$ by the natural extension of the action 
$t^m \cdot \zeta_{\f} = \zeta_{\overline{t^mf}}$.
We note that the whole of $L$ cannot act on $M$ in this way, since we would get
contradictions such as 
$0 = 2p \cdot \zeta_p = 2 \zeta_{p^2} = 2^{r(p^2)} \zeta_{p^2} = p^2 \neq 0$.
In particular, we do not have a ring structure on $M$.
For each $n\geq 0$ we denote by $K_n \su L_n \su M_n$ the additive subgroups
of $K\su L \su M$ respectively generated by the monomials of degree $n$
(where $\zeta_\f \in M$ is considered a monomial of the same degree as $\f$).
We note $L_1=M_1=\G_U$.
We have $L_1=K_1\oplus S$ where $S \su L_1$
is the four element subgroup generated by $p,q$.

We now define a function
$\FF:L_1\to M$ as follows: We first define $\FF:K_1 \to K$ as the group 
homomorphism
from the additive group $K_1 = \{ mt : m \in \Z \} \cong \Z$ to the multiplicative group of 
invertible elements in $K$, given on the generator $t$ of $K_1$ 
by $$\FF(t)=\sum_{n=0}^\infty t^n.$$ This is indeed an invertible element, giving
$\FF(-t)=\bigg(\sum_{n=0}^\infty t^n\bigg)^{-1} = 1-t$.

We then define $\FF:S\to M$ explicitly on the four elements of $S$ as follows: 
\begin{enumerate}
\item $\FF(0)=1$.
\item $\FF(p)=\sum_{n=0}^\infty \zeta_{p^n}$.
\item $\FF(q)=\sum_{n=0}^\infty \zeta_{q^n}$.
\item $\FF(p+q)=1+p+q+\sum_{n=2}^\infty (\zeta_{p^n} + \zeta_{q^n} + \zeta_{\overline{pq^{n-1}}})$.
\end{enumerate}

Finally, $\FF:L_1 \to M$ is defined as follows: Any element in $L_1$ is uniquely written as
$k+s$ with $k\in K_1, s\in S$, and we define $\FF(k+s)=\FF(k)\FF(s)$
where the product on the right is the action of $K$ on $M$. 
Let $\FF_n:L_1 \to M_n$ be the projection of $\FF$ into $M_n$.
We now state our classification theorem for finite order invariants.
The proof follows exactly as for the corresponding
claim for orientable surfaces appearing in \ct{h}. 

\begin{thm}\label{high}
For any closed non-orientable surface $F$, regular homotopy class $\A$ of immersions of $F$ 
into $\E$ and Abelian group $\G$, the 
image of the injection $u:V_n / V_{n-1} \to \Delta_n$ is $E_n$,
and the invariant $\FF_n \circ f_1^U : I_0 \to M_n$ is a universal order $n$ invariant.
\end{thm}

\section{$\4$-forms and automorphisms of $F$}\label{og}

Denote $\4=(\hf\Z)/(2\Z)$, which is a cyclic group of order 4 
(and $\4$ stands for \emph{half} integers).
The group $\C=\Z / 2\Z$ is contained in $\4$ as a subgroup.

\begin{dfn}\label{d1}
Let $E$ be a finite dimensional vector space over $\C$.
An $\4$-form is a map $g:E \to \4$ satisfying the following two conditions:

\begin{enumerate}
\item  For all $x,y \in E$: $$g(x+y) = g(x) + g(y) + C(x,y)$$
where $C:E \times E \to \C$  ($\su\4)$ is a non-degenerate symmetric bilinear form.
\item There is at least one $x \in E$ with $g(x) \not\in \C$
i.e. $g(x)=\pm\hf$.
\end{enumerate}
\end{dfn}

It follows that $g(0)=0$ and $2g(x) = C(x,x)$ showing $C(x,x) =1$ iff
$g(x) = \pm\hf$ and so there is at least one $x \in E$ with $C(x,x)=1$.
One can then show
there exists an ``orthonormal'' basis for $E$, i.e. a basis $e_1,\dots,e_n$ satisfying
$C(e_i,e_j)=\delta_{ij}$.
For such basis if $x=\sum_{j=1}^k e_{i_j}$ ($1\leq i_1 < \cdots < i_k \leq n$) 
then $g(x)=\sum_{j=1}^k d_{i_j}$ where $d_i=g(e_i)=\pm\hf$.

For $\4$-form $g$ on $E$ we define $O(E,g)$ to be the group of all
linear maps $T:E \to E$ satisfying $g(Tx)=g(x)$ for all $x \in E$.
It follows that $C(Tx , Ty) = C(x, y)$ for all $x,y \in E$ and that 
$T$ is invertible.
For $a \in E$ we define $T_a:E \to E$ to be the map $T_a(x)=x + C(x,a) a$.
Then $T_a \in O(E,g)$ iff $g(a)=1$ or $a=0$. 
For $a,b \in E$ we define $S_{a,b}:E \to E$ by:
$$S_{a,b}(x) =   T_a \circ T_b \circ T_{a+b} .$$  
One verifies directly that
if $a,b \in E$ satisfy $g(a)=g(b)=g(a+b)=0$ (which is equivalent to
$g(a)=g(b)=C(a,b)=0$), then $S_{a,b} \in O(E,g)$. 
In this case $S_{a,b}$ may also be written as: $S_{a,b}(x)=x + C(x,b)a + C(x,a)b$.

\begin{thm}\label{orth}
Let $g$ be an $\4$-form on $E$, then $O(E,g)$ is generated by
the set of elements of the following two forms:
\begin{enumerate}
\item $T_a$ for $a\in E$ with $g(a)=1$.
\item $S_{a,b}$ for $a,b \in E$ with $g(a)=g(b)=g(a+b) = 0$.
\end{enumerate}
Furthermore, if $\dim E\geq 9$
then the elements of the first form alone generate $O(E,g)$.

\end{thm}

\begin{pf}
Let $T\in O(E,g)$ and let $e_1,\dots,e_n$
be an orthonormal basis for $E$. Assuming inductively, with decreasing $k$,
that $T$ fixes $e_{k+1},\dots,e_n$, we will compose 
$T$ with elements of the above two forms to obtain a map which additionally
fixes $e_k$, eventually fixing all $e_i$ which will prove the statement.
For $x \in E$ define $\spp(x)$ to be the
set $\{ e_{i_1},\dots,e_{i_m} \}$ 
($i_1 < \cdots < i_m$) such that $x=\sum_{j=1}^m e_{i_j}$. 
For $i \leq k$ and $j \geq k+1$, $C(Te_i,e_j)=C(e_i,e_j)=0$, 
and so for $i \leq k$, $\spp(Te_i) \su \{e_1,\dots,e_k \}$.
Denote $v=Te_k$, then $g(v)=g(e_k)=\pm\hf$.

\emph{Case A}: $C(v, e_k) =0$.
Letting $a=e_k+v$ we have $g(a)=g(e_k)+g(v)+C(e_k, v) =1$, 
$T_a \circ T(e_k)=T_a(v)=e_k$, and for $i \geq k+1$, $T_a \circ T (e_i)=T_a(e_i)=e_i$, so we are done.

\emph{Case B}: $C(v, e_k) =1$.
Since $\spp(v) \su \{e_1,\dots,e_k \}$, if $k=1$ then $v=e_k$ and we are done, and 
so we assume $k \geq 2$. 
We first show that $\spp(v) \neq \{e_1,\dots,e_k \}$.
Indeed if $v=\sum_{i=1}^k e_i$ then 
$C(Te_{k-1} ,\sum_{i=1}^k e_i) = C(Te_{k-1} , Te_k) =C(e_{k-1} , e_k)=0$.
Since $\spp(Te_{k-1}) \su \{ e_1,\dots,e_k \}$ this shows
that the number of elements in $\spp(Te_{k-1})$
is even, which implies $C(Te_{k-1} , Te_{k-1})=0$,
contradicting $C(e_{k-1} , e_{k-1}) =1$. 
Since $C(v , e_k) =1$ we know $e_k \in \spp(v)$, so say $e_1 \not\in \spp(v)$. 
If $g(e_1)=g(e_k)$ then 
since $C(e_1, v) =0$ and $C(e_1 , e_k) =0$,
by the argument of Case A, we can map $v$ to $e_1$ and then $e_1$ to $e_k$ and we are done. 
Otherwise, say $g(v)=\hf$ and $g(e_1)=-\hf$. 
Assuming $v \neq e_k$ (otherwise we are done), $e_k$ is not the only element in
$\spp(v)$.
The condition $g(v)=g(e_k)$ implies one of three possibilities (after relabeling indices): 
\begin{enumerate}
\item $k\geq 4$, $e_2,e_3 \in \spp(v)$, $g(e_2)=\hf$, $g(e_3)=-\hf$.
\item $k\geq 6$, $e_2,e_3,e_4,e_5  \in \spp(v)$ and the value of $g$ on
$e_2,e_3,e_4,e_5$ is $\hf$.  
\item $k\geq 6$, $e_2,e_3,e_4,e_5  \in \spp(v)$ and the value of $g$ on
$e_2,e_3,e_4,e_5$ is $-\hf$.
\end{enumerate}
In case 1 define $a=e_1+e_2, b=e_3+e_k$.
In case 2 define $a=e_1+e_2, b=e_3+e_4+e_5+e_k$.
In case 3 define $a=e_1+e_2+e_3+e_4, b=e_5+e_k$.
In all three cases $g(a)=g(b)=g(a+b)=0$ and so $S_{a,b}$ belongs to the
set of proposed generators.
In all three cases $C(v , a) =1$ and $C(v, b) =0$ and so 
$S_{a,b}(v)=v+b$. Furthermore, $S_{a,b}$ fixes all $e_j$ for $j \geq k+1$, 
and $C(v+b , e_k) = C(v, e_k) + C(b, e_k) = 1 + 1 =0$.
So we can use $S_{a,b}$ to map $v$ to $v+b$, 
and then by Case A we can map $v+b$ to $e_k$, and we are done.

We conclude by showing that if $n \geq 9$ then the elements $T_a$ alone generate
$O(E,g)$. Indeed in the proof above we used only maps $S_{a,b}$ where  
$\spp(a) \cup \spp(b)$
includes at most six elements. If $n \geq 9$ then for
each such pair $a,b$ there are at
least three basis elements not in $\spp(a) \cup \spp(b)$. 
In the span of such three basis elements there is an element $s$ with $g(s)=1$. 
One verifies directly, by checking on basis elements, 
that $S_{a,b} = T_s \circ T_{s+a} \circ T_{s+b} \circ T_{s+a+b}$.
Since $g(s)=g(s+a)=g(s+b)=g(s+a+b)=1$ we have indeed expressed the given element
$S_{a,b}$ as a product of four generators of the first form.
\end{pf}

\begin{dfn}\label{d2}
For $T \in O(E,g)$ let $\psi(T) = \psi_E(T)= \r(T - \I) \bmod{2} \in \C$
\end{dfn}

\begin{prop}\label{hom}
The map $\psi : O(E,g) \to \C$ is a homomorphism. 
\end{prop}

\begin{pf}

Assume first that $\dim E \geq 9$. In this case we know $O(E,g)$ is generated by elements 
of the form $T_a$ with $g(a)=1$. 
Let $\F(T) = \{x \in E : Tx=x\} = \ker(T-\I)$ then $\r(T-\I) = \cd \F(T)$. 
We may now use \ct{a} Lemma 3.1, which is stated in a slightly different setting, 
(of $\C$ valued quadratic forms,) but whose proof applies word by word to our setting.
The statement of \ct{a} Lemma 3.1 is as follows: If $g(a)=1$ then for any $T\in O(E,g)$,
$\cd \F(T \circ T_a) = \cd \F(T) \pm 1$. Since $\cd \F(T_a) =1$ and the
elements $T_a$ generate $O(E,g)$, $\psi$ is a homomorphism.

Assume now $\dim E <9$ and take some $E'$ with $\4$-form $g'$ such that
$\dim E + \dim E' \geq 9$. The function
$g\oplus g' : E \oplus E' \to \4$ defined by 
$g \oplus g' (x,x') = g(x)+g'(x')$, is an $\4$-form on $E \oplus E'$.
Let $u:O(E,g) \to O(E\oplus E' , g \oplus g')$ be the embedding given by
$u(T)= T \oplus \I_{E'}$. Since $\r(T-\I_E) = \r(T \oplus \I_{E'} - \I_{E \oplus E'})$
we have $\psi_E = \psi_{E \oplus E'} \circ u$.
Since $\dim(E \oplus E') \geq 9$,
$\psi_{E \oplus E'}$ is a homomorphism and so $\psi_E$ is a homomorphism.
\end{pf}

Returning to surfaces,
let $F$ be a closed non-orientable surface, and let $g:\hc \to \4$ be an $\4$-form
whose associated bilinear form $C(x,y)$ is the algebraic intersection form $x \cdot y$
on $\hc$. Let $\N = \N(F)$ be the group of all diffeomorphisms $h:F \to F$ up to isotopy.
For $h:F \to F$ let $h_*$ denote the map it induces on $\hc$. The subgroup
$\ng = \N(F)_g \su \N$ is defined by
$\ng = \{ h \in \N : h_* \in O(\hc , g) \}$.

A simple closed curve will be called a \emph{circle}. If $c$ is a circle in $F$, 
the homology class of $c$ in $\hc$ will be denoted by $[c]$. A circle $c$ in $F$ 
has an annulus neighborhood if $[c] \cdot [c] =0$ and a Mobius band neighborhood 
if $[c]\cdot[c] =1$. Such circles will be called $A$-circles and $M$-circles, respectively.
Given an $A$-circle $c$ in $F$, a Dehn twist along $c$ will be denoted $\Tc$.
The map induced on $\hc$ by $\Tc$ is $T_{[c]}$,
and so $\Tc \in \ng$ iff $g([c])=1$ or $[c]=0$.
Also, since $(T_{[c]})^2=\I$ whenever $[c]\cdot [c] =0$, $(\Tc)^2 \in\ng$ for any $A$-circle $c$.

Let $P \su F$ be a disc with two holes. 
Let $c,d,e$ be the three boundary circle of $P$, then $[c]+[d]=[e]$.
If $g([c])=g([d])=0$ (and so $g([e])=0$) then define $\SP = \Tc \circ \T_d \circ \T_e$.
The map induced by $\SP$ on $\hc$ is $S_{[c],[d]}$.
Finally, a $Y$-map $h:F\to F$ as defined in \ct{l} induces the identity on
$\hc$ and so $h \in \ng$.

\begin{dfn}\label{d3}
Let $F$ be a closed non-orientable surface and $g$ an $\4$ form on $\hc$. 
A map $h:F \to F$ will be called \emph{good} if it is of one of the following five
forms:
\begin{enumerate}
\item $h = (\Tc)^2$ for any $A$-circle $c$. 
\item $h=\Tc$ for an $A$-circle $c$ with $g([c])=1$
\item $h=\Tc$ for an $A$-circle $c$ with $[c]=0$
\item $h = \SP$ for some disc with two holes $P \su F$ with boundary circles
$c,d,e$ satisfying $g([c])=g([d])=0$ ($=g([e])$.
\item $h$ is a $Y$-map (defined in \ct{l}).
\end{enumerate}
A good map will be called of type 1 - 5 accordingly.

\end{dfn}

Whenever we consider two circles in $F$, we will assume they intersect transversally,
$|c \cap d|$ will then denote the number of intersection points between circles $c,d$.
(And so the algebraic intersection $[c]\cdot [d]$ in $\hc$ is the reduction mod 2 of 
$|c \cap d|$.)

\begin{thm}\label{K}
Let $F$ be a closed non-orientable surface and let $K_F = \{ h \in \N : h_* = \I \}$.
Then $K_F$ is generated by the good maps of type 1,3,5.
\end{thm}

\begin{pf}
Let $c_1,\dots,c_n$ be a family of disjoint $M$-circles in $F$ such that cutting
$F$ along $c_1,\dots,c_n$ produces a disc with $n$ holes. Then $[c_1],\dots,[c_n]$
is an orthonormal basis for $\hc$.
Now let $h \in K_F$.
We will compose $h$ with maps of type 1,3,5 until we obtain the identity.
Assume inductively that $h(c_i)=c_i$ for all $i \leq k-1$ (not necessarily respecting
the orientation of $c_i$). Denote $a = h(c_k)$, then $a$ is disjoint from $c_1,\dots,c_{k-1}$
and since $[a] \cdot [c_j] = [c_k] \cdot [c_j] = 0$ for $j \geq k+1$,
$|a \cap c_j|$ is even for $j \geq k+1$. We will now perform good maps of type 1,
which fix $c_1,\dots,c_{k-1}$, to map $a$ to a circle disjoint from all $c_j$, $j \geq k+1$.
Assume for some $j \geq k+1$, $|a \cap c_j|>0$, then since it is even and since $c_j$ is 
an $M$-circle, there must be two adjacent intersection
points along $c_j$ at which $a$ crosses $c_j$ in the same direction with respect to an
orientation on a neighborhood of the subinterval $b \su c_j$ satisfying $b \cap a = \pa b$.
Now $\pa b$ divides $a$ into two intervals $a',a''$ and since $a$ is an $M$ circle,
one of the circles $a' \cup b$ or $a'' \cup b$ is an $A$-circle. Say $c = a' \cup b$ is 
an $A$-circle. Then as shown in \ct{a} Figure 2, 
$|(\Tc)^2(a) \cap c_j|=|a \cap c_j|-2$. 
This map fixes 
$c_1,\dots,c_{k-1}$ and we may repeat this procedure until the image of $a$ (which we again
name $a$), is disjoint from all $c_j$, $j \neq k$.

We will now use good maps of type 5 to map $a$ onto $c_k$. 
Let $M_1,\dots,M_{k-1},M_{k+1},\dots,M_n$ be disjoint Mobius band neighborhoods 
of $c_1,\dots,c_{k-1},c_{k+1},\dots,c_n$ which are also disjoint from $c_k$ and $a$. 
Let $\tF$ be the projective plane obtained from $F$ by collapsing each $M_i$ to a point
$p_i$ ($i=1,\dots,k-1,k+1,\dots,n$). In $\tF$ one can isotope $a$ to coincide with $c_k$.
We lift this isotopy to $F$ as follows: Every time $a$ is about to pass one of the points $p_i$
in $\tF$, we realize this passage in $F$ by performing a good map of type 5, which drags
$M_i$ along a circle $a'$ which is close to $a$ and intersects $a$ once, returning $M_i$ 
to place (while reversing the orientation of $c_i$). This pushes $a$ to the other side 
of $\pa M_i$ in $F$ and so to the other side of $p_i$ in $\tF$.
In this way we bring $a$ to coincide with $c_k$.

We continue by induction until $h(c_i)=c_i$ for all $1\leq i \leq n$, but when performing the 
procedure involving $Y$-maps with the last circle $c_n$, we choose the isotopy in $\tF$
to bring $h(c_n)$ \emph{orientation preservingly} onto $c_n$. Since cutting $F$ along
$c_1,\dots,c_n$ produces an orientable surface (a sphere with $n$ holes), then if $c_n$ 
is mapped onto itself orientation preservingly, then the same must be true for
all $c_i$. After some isotopy we obtain that $h$ is the identity on all $c_i$.

Let $S$ be the sphere with $n$ holes obtained from $F$ by cutting it along 
$c_1,\dots,c_n$, then $h$ induces a diffeomorphism from 
$S$ to $S$ which is the identity on $\pa S$.
It is known that any such map is a composition of Dehn twists. Since any circle 
in $S$ separates $S$, the corresponding circle in $F$ will separate $F$ and so these
Dehn twists on $F$ are good maps of type 3.

\end{pf}

\begin{thm}\label{ng}
Let $F$ be a closed non-orientable surface and let $g$ be an $\4$-form on $\hc$.
Then $\ng$ is generated by the good maps. 
Furthermore, if $\dim \hc \geq 9$ then good maps of type 4 are not needed.
\end{thm}

\begin{pf}
For $h \in \ng$ we first compose $h$ with good maps of type 2,4 to obtain a map
in $K_F$, as follows. Any element $a \in \hc$ with $g(a)=1$ may be realized by
a circle in $F$ and so any generator $T_a$ of Theorem \pr{orth} may be realized 
by a good map of type 2.
Any $a,b \in \hc$ with $g(a)=g(b)=a \cdot b =0$ may be realized by a pair of disjoint circles
$c,d$ in $F$. Let $l$ be an arc connecting $c$ to $d$ (its interior being disjoint from $c,d$)
then a regular neighborhood $P$ of $c \cup l \cup d$ is a disc with two holes with $c,d$
isotopic to two of its boundary components. So any generator $S_{a,b}$ of Theorem \pr{orth}
may be realized by a good map of type 4. So by Theorem \pr{orth} we may compose $h$ with 
good maps
of type 2,4 to obtain a map in $K_F$, and if $\dim \hc \geq 9$ then maps of type 4 
are not needed. We complete the proof using Theorem \pr{K}.
\end{pf}

\section{Tangencies and quadruple points of regular homotopies}\label{fmq}

For regularly homotopic immersions $i,j:F \to \E$, let $P(i,j)=P(i)-P(j)$ and $Q(i,j)=Q(i)-Q(j)$.
So $P(i,j) \in \C$ (respectively $Q(i,j) \in\C$)
is the number mod 2 of tangency points (respectively quadruple points)
occurring in any generic regular homotopy between $i$ and $j$. 
In this section we prove the following:

\begin{thm}\label{qm}
Let $F$ be a closed non-orientable surface. Let $i:F\to \E$ be a stable immersion
and let $h:F \to F$ be a diffeomorphism such that $i$ and $i \circ h$ are regularly
homotopic.
Then 
$$P(i,i\circ h)=Q(i,i\circ h) =\bigg(\r(h_*-\I)+\ep(\det h_{**})\bigg)\bmod{2}$$
where $h_*$ is the map induced by $h$ on $\hc$, $h_{**}$ is
the map induced by $h$ on $H_1(F;\Q)$ and for $0 \neq q \in \Q$, 
$\ep(q) \in \C$ is 0 or 1
according to whether $q$ is positive or negative, respectively.
\end{thm}

For $h \in \ng$ define 
$\Omega(h)= ( \r(h_*-\I)+\ep(\det h_{**}) )\bmod{2} =\psi(h_*) + \ep(\det h_{**})  \in \C$.
So we must show $P(i,i\circ h)=Q(i,i\circ h)=\Omega(h)$.
By \ct{f} Proposition 7.1 we know $P(i,i\circ h)=Q(i,i\circ h)$. 
Indeed the proof of this case appearing there, 
does not use any assumption on orientability of the surface.
So it is enough to prove $Q(i,i\circ h)=\Omega(h)$.

Let $F$ be any closed surface and let $c \su F$ be a separating circle. Denote by
$F_1,F_2$ the two subsurfaces into which $c$ divides $F$. Denote by $\one$ the closed surface
obtained from $F$ by collapsing $F_2$ to a point, and similarly denote by $\two$
the closed surface obtained from $F$ by collapsing $F_1$ to a point. 
If $h:F \to F$ is a diffeomorphism such that $h(F_1) = F_1$ and $h(F_2) = F_2$
then $h$ induces maps $h_k:\tfk \to \tfk$ ($k=1,2$). 
As appears in \ct{a} Section 5.3, 
if $i:F \to \E$ is an immersion then it determines 
regular homotopy classes of immersions  $i_k:\tfk \to \E$, 
and if $h$ as above satisfies that $i$ and $i \circ h$ are regularly homotopic
then also $i_k$ and $i_k \circ h_k$ are regularly homotopic ($k=1,2$).
It is shown in \ct{a} Section 5.3 for orientable $F$, that for such $i$ and $h$,
if $h$ is orientation preserving then 
$Q(i,i\circ h)=Q(i_1,i_1\circ h_1) + Q(i_2,i_2\circ h_2)$ ($\in \C$)
whereas if $h$ is orientation reversing then
$Q(i,i\circ h)=Q(i_1,i_1\circ h_1) + Q(i_2,i_2\circ h_2) +1$ ($\in \C$).
But in fact the proof does not use the property that $h$ itself is orientation preserving
or reversing, but rather the corresponding property that $h|_c:c \to c$ is orientation
preserving or reversing. This latter property is meaningful also for non-orientable
surfaces, and indeed the proof works just the same for non-orientable surfaces. And so we have:

\begin{lemma}\label{add}
Let $F$ be any closed surface. In the above setting, if
$h|_c:c\to c$ is orientation preserving then 
$$Q(i,i\circ h)=Q(i_1,i_1\circ h_1) + Q(i_2,i_2\circ h_2)$$
and if $h|_c:c\to c$ is orientation reversing then
$$Q(i,i\circ h)=Q(i_1,i_1\circ h_1) + Q(i_2,i_2\circ h_2) +1.$$
\end{lemma}
 
\begin{cor}\label{cadd}
In the above setting if $h$ satisfies that $h(x)=x$ for all $x \in F_2$
then $$Q(i,i\circ h)=Q(i_1,i_1\circ h_1).$$ 
\end{cor}

\begin{lemma}\label{tf}
In the above setting if $h$ satisfies that $h(x)=x$ for all $x \in F_2$,
then $$\Omega(h)=\Omega(h_1).$$ 
\end{lemma}

\begin{pf}

We first show $\psi(h_*)=\psi({h_1}_*)$. We have 
$\hc \cong H_1(F_1;\C) \oplus H_1(F_2;\C) \cong H_1(\one;\C) \oplus H_1(\two;\C)$ and
under this isomorphism, $h_*$ corresponds to ${h_1}_* \oplus \I$. It follows
that $\psi(h_*)=\psi({h_1}_*)$. 

Next we show $\det h_{**} = \det {h_1}_{**}$. Since
$\widetilde{H}_0(F_2;\C)=0$, we have the 
exact sequence 
$$H_1(F_2;\Q) \to H_1(F;\Q) \to H_1(F,F_2;\Q) \to 0.$$
Let $r$ denote the map that $h$ induces on $H_1(F,F_2;\Q)$, then since
$h$ induces the identity on $H_1(F_2;\Q)$, we have $\det h_{**} = \det r$.
But $r$ corresponds to ${h_1}_{**}$ under the natural isomorphism
$H_1(F,F_2;\Q) \to H_1(\one,p;\Q) = H_1(\one;\Q)$.

\end{pf}

Now for non-orientable $F$
let $i:F \to \E$ be an immersion and let $g^i$ be the $\4$-form on $\hc$ 
determined by $i$, as defined in \ct{p}. 
Note that the notation in \ct{p} differs from ours in
that $\4$ is taken there to be $\Z / 4\Z$ rather than $(\hf\Z) / 2\Z$. And so
the numerical value of $g^i$ appearing there is twice the value here, and
our relation $g(x+y)=g(x)+g(y)+x \cdot y$ is replaced there by 
$g(x+y)=g(x)+g(y)+ 2(x \cdot y)$.
We have as in \ct{a} Proposition 5.2, for any diffeomorphism 
$h:F \to F$, $i$ and $i \circ h$ are regularly homotopic iff $h \in \N_{g^i}$.

Letting $n=\dim \hc$, we
will prove that $Q(i,i \circ h) = \Omega(h)$ for all $n$ in the following order:
$n=1$, $n=2$, $n \geq 9$, and finally $3 \leq n \leq 8$. 
Denote $g=g^i$.
Since by Theorem \pr{hom}, $\psi$ is a homomorphism 
on $O(\hc,g)$
and so $\Omega$ is a homomorphism on $\ng$, and since by \ct{a} Lemma 5.5
(whose proof applies to non-orientable $F$) $h \mapsto Q(i,i \circ h)$ is a homomorphism,
it is enough to check the equality $Q(i,i \circ h) = \Omega(h)$ for generators
of $\ng$. 

For $n=1$, $F=\R P^2$ the projective plane. Since $\N(\R P^2)$ is trivial,
the equality follows trivially.

For $n=2$, $F=Kl$ the Klein bottle. It is shown in \ct{l} that
$\N(Kl) \cong \C \oplus \C$. There is a unique circle $c$ in $Kl$ up to isotopy which
separates $Kl$ into two punctured projective planes $F_1,F_2$. 
We have $\dim H_1(Kl;\C) =2$ with unique orthonormal basis $e_1,e_2$ 
where $e_k$ ($k=1,2$) corresponds to the unique generator of $H_1(F_k,\C)$. 
With respect to this basis 
$h_* \in \{ 
\left(\begin{smallmatrix} 1 & 0 \\ 0 & 1 \end{smallmatrix}\right) ,
\left(\begin{smallmatrix} 0 & 1 \\ 1 & 0 \end{smallmatrix}\right) \}$
for all $h \in \N(Kl)$.
On the other hand $\dim H_1(Kl;\Q) =1$ with generator given by the circle $c$, and
we have $h_{**} \in \{ \I , -\I \}$ for all $h \in \N(Kl)$. 
The map $h \mapsto (h_*,h_{**})$ indeed realizes the 
isomorphism $\N(Kl) \cong \C \oplus \C$. 
For the unique $Y$-map $u:Kl \to Kl$, we have $u_* = \I$ and so 
$u \in \ng$ and $\psi(u_*)=0$. On the other hand $u_{**}=-\I$ so 
$\det u_{**} =-1$ so $\ep(\det u_{**}) =1$ giving finally $\Omega(u)=1$.
Now we may isotope $u$ such that $u(F_1)=F_1$ and $u(F_2)=F_2$ and
in this case $u|_c : c \to c$ is orientation reversing. 
Since for $k=1,2$, $u_k \in \N(\R P^2)$ which is trivial, 
we have by Lemma \pr{add} that $Q(i,i\circ u)=1$. 
So we have shown $Q(i,i\circ u)=1=\Omega(u)$.

If $g(e_1) \neq g(e_2)$ then $u$ is the only non-trivial map in $\N_g$ and so we are done.
Otherwise $g(e_1)=g(e_2)$ and then $\ng=\N$. By \ct{a} Lemma 5.8
(whose proof applies to non-orientable $F$) we may replace $i$ by any other immersion in
its regular homotopy class. Indeed we construct $i$ as follows: $F_1$ will be immersed
in the positive side of the $yz$ plane, with its boundary $c$ embedded in the $yz$ plane
symmetrically with respect to the $z$ axis. The immersion of $F_1$ is furthermore chosen to give
the correct value for $g(e_1)$. Now the image of $F_2$ will be the punctured projective 
plane obtained from the image of $F_1$ by a $\pi$ rotation around the $z$ axis, and
so $g(e_2)=g(e_1)$ also has the correct value and so the new immersion $i$ is indeed
in the correct regular homotopy class.
The symmetry of $i(F)$ implies that there is
a diffeomorphism $v:F \to F$ such that $\pi$ rotation of $i(F)$ around the $z$ axis is a 
regular homotopy between $i$ and $i \circ v$. 
Now $v_* = \left(\begin{smallmatrix} 0 & 1 \\ 1 & 0 \end{smallmatrix}\right)$
and so $v$ together with the $Y$-map $u$ generate $\N$, and $\psi(v_*)=1$.
Since $v$ reverses the orientation of $c$ we have $v_{**}=-\I$ and so
$\ep(\det v_{**})=1$.
Together, $\Omega(v)=1+1=0$. Now the rigid rotation of $i(F)$ 
which is the regular homotopy between $i$ and $i \circ v$, has no quadruple points at all, and 
so we get $Q(i,i\circ v)=0=\Omega(v)$, which completes the case $n=2$.

For $n\geq 9$ we have $\ng$ generated by good maps of type 1,2,3,5.
If $h \in \ng$ is of type 3, then $h=\Tc$ for a separating circle $c \su F$.
Let $c'$ be a circle parallel to $c$, dividing $F$ into $F_1,F_2$ with $c \su F_1$.
Then $h(x)=x$ for all $x \in F_2$ and $h_1$ is isotopic to the identity
on $\one$. And so by Corollary \pr{cadd} and Lemma \pr{tf} we have
$Q(i,i\circ h) = Q(i_1,i_1\circ h_1) = 0 = \Omega(h_1) = \Omega(h)$.
If $h$ is of type 1,2, then $h$ is $\Tc$ or $(\Tc)^2$ where (since
we have already discussed good maps of type 3), we may assume $c$
is a non-separating $A$-circle. So there is an $M$-circle $d \su F$ such that
$|c \cap d|=1$ and so a regular neighborhood $F_1$ of $c \cup d$ is a punctured
Klein bottle. Since we have already established the theorem for $n=2$ we know
$Q(i_1,i_1\circ h_1)= \Omega(h_1)$. Since $h(x)=x$ for all $x \in F_2 = F-F_1$,
we have by Corollary \pr{cadd} and Lemma \pr{tf}, 
$Q(i,i\circ h) = Q(i_1,i_1\circ h_1) = \Omega(h_1) = \Omega(h)$.
Finally if $h$ is of type 5, i.e. a $Y$-map, then again there is a punctured
Klein bottle $F_1 \su F$ with $h(x)=x$ for $x \in F-F_1$ and so we are done as in the
previous case, which completes the proof for $n\geq 9$.

For $3 \leq n \leq 8$, since $F$ is non-orientable, $h$ is isotopic to a map
which fixes a disc $D \su F$ pointwise, so assume $h$ satisfies this property.
Take any closed surface $F'$ with 
$\dim \hc + \dim H_1(F';\C) \geq 9$ and construct $G = F \# F'$ where the connect
sum operation is performed by deleting the above mentioned disc $D$ from $F$, and some
disc $D'$ from $F'$. There clearly exists an immersion $j:G \to \E$ such that 
$j|_{F-D}=i|_{F-D}$. 
By our assumption on $D$ we can extend $h|_{F-D}$ to a diffeomorphism
$u:G \to G$ by defining $u(x)=x$ for all $x \in F' -D'$. 
Denote $G_1=F-D$ and $G_2 = F'-D'$ then $\widetilde{G_1}$ is naturally identified with
$F$ and under this identification $j_1$ corresponds to $i$ and $u_1$ corresponds to $h$.
Since $\dim H_1(G;\C) \geq 9$ our theorem is already proved for $G$ and so by Corollary
\pr{cadd} and Lemma \pr{tf} we get
$Q(i,i\circ h) = Q(j,j\circ u) = \Omega(u) = \Omega(h)$.
This completes the proof of Theorem \pr{qm}.

\end{document}